\documentclass[12pts]{amsart}
\usepackage{amssymb, amsmath}
\usepackage{amsthm}
\usepackage{amscd}
\usepackage{graphicx}
\usepackage[ngerman]{babel}

\newtheorem{theorem}{Satz}

\theoremstyle{definition}

\newtheorem*{remark}{Bemerkung}

\title{Eine Bemerkung zu einigen $2$-dimensionalen Komplexen, die im $\mathbb{R}^4$ fast-eingebettet werden k\"onnen}
\author{T. T$\hat{\mathrm{a}}$m Nguy$\tilde{\hat{\mathrm{e}}}$n-Phan}
\address{Karlsruher Institut f\"ur Technologie\\
Karlsruhe\\
Deutschland}
\email{tu.phan@kit.edu}

\def\beqa{\begin{eqnarray}}
\def\eeqa{\end{eqnarray}}
\def\beqa{\begin{eqnarray}}
\def\eeqa{\end{eqnarray}}

\input xy
\xyoption{all}

\begin{document}

 
\maketitle

\large{
Auf die Frage, wann ein gegebener $2$-dimensionaler Komplex im $\mathbb{R}^4$ eingebettet werden kann, gibt es bisher keine vollst\"andige Antwort. Nach der Arbeit von van Kampen  \"uber die allgemeine Frage der Einbettungsm\"oglichkeit eines $n$-dimensionalen Komplexes $X^n$ im $\mathbb{R}^{2n}$, k\"onnen wir nur bestimmen, ob es eine abteilungsweise lineare Immersion  $f\colon X^n \rightarrow \mathbb{R}^n$
gibt, so dass die algebraische Schnittpunktszahl jedes Paares getrennter\footnote{Zwei $n$-dimensionale Simplices von $X^n$ hei"sen \emph{getrennt}, wenn sie keinen gemeinsamen Punkt in $X$ haben.} Simplices von $X^n$ gleich Null ist (\cite{vanKampenpaper}, \cite{vanKampenerr}). In einem solchen Fall sagen wir, dass die \emph{van-Kampen-Obstruktion} f\"ur $X$ verschwindet, und wenn $n$ auch gr\"o"ser als $2$ ist, k\"onnen wir mit Hilfe einiger Whitney-Disketten alle Singularit\"aten von $f$ entfernen, um eine Einbettung $X^n \rightarrow \mathbb{R}^{2n}$ zu erhalten. Aber wenn $n=2$, ist es nicht immer m\"oglich Whitney-Disketten in $\mathbb{R}^4$ zu finden, und tats\"achlich kann man nicht immer alle Singularit\"aten von $f$ entfernen, wie die Beispiele $2$-dimensionaler Komplexe von Freedman-Kruskal-Teichner, die im $\mathbb{R}^4$ nicht eingebettet werden k\"onnen, aber deren van-Kampen-Obstruktionen verschwinden, zeigen.} 

\subsection*{Die Freedman-Kruskal-Teichner-Beispiele (\cite{FKT})} \large{
Sei $X$ ein $2$-dimensionaler Komplex aufgebaut aus allen $2$-dimensionalen Simplices, deren Eckpunkte zu $\{x_0, x_1, ..., x_6\}$ geh\"oren, mit Ausnahme derer mit Eckpunkten  $x_4, x_5, x_6$. Sei $Y$ der gleiche Komplex wie $X$, aber mit Eckpunkten in $\{y_0,  ..., y_6\}$. Sei $Z$ der $2$-Komplex, den wir aus der disjunkten Vereinigung $X  \sqcup Y$ durch die Identifizierung der Punkte $x_6 \sim y_6$ erhalten, deren Bild in $Z$ wir als Basispunkt $o$ von $Z$ w\"ahlen. Wir bezeichnen durch $a$ (bzw. $b$) die geschlossene Kurve in $Z$ mit Basispunkt $o$, die aus den drei Kanten $ox_4$, $x_4x_5, x_5o$ (bzw. $oy_4$, $y_4y_5, y_5o$) aufgebaut wird. Schlie"slich bezeichnen wir durch $F = \langle a, b\rangle$ die freie Gruppe in $\{a,b\}$. Die Komplexe von Freedman-Kruskal-Teichner sind die Komplexe      
$$ K_\varphi = Z \cup D^2,$$
die durch Anheften einer $2$-dimensionalen Zelle $D^2$ an $Z$ mittels der Anheftungsabbildung $\varphi$ entstehen, so dass $\varphi \colon \partial D^2 \rightarrow a\vee b$  in der Kommutatorgruppe $[F,F]$ liegt und $\varphi \ne 1$.
\newline}

\large{
Es ist einfach, eine Immersion $K_\varphi \rightarrow \mathbb{R}^4$ zu konstruieren, deren einzige Singularit\"aten  Schnittpunkte der Zellen $x_1 x_2 x_3$ und $y_1 y_2 y_3$ sind, deren algebraische Schnittpunktszahl gleich Null ist. Viele dieser Komplexe $K_\varphi$ k\"onnen tats\"achlich in $\mathbb{R}^4$ \emph{fast-eingebettet} werden, d.h. es gibt eine (abteilungsweise lineare) Abbildung  $g\colon K_\varphi \rightarrow \mathbb{R}^4$, so dass die Bilder getrennter Simplices einander nicht schneiden. Tats\"achlich beweisen wir den folgenden st\"arkeren Satz.}

\begin{theorem}
\large{Wenn $\varphi \in F^{(3)} = [F,[F,F]]$, kann der Komplex $K_\varphi$ im $\mathbb{R}^4$ immersiert werden, so dass die einzigen Singularit\"aten Selbstschnittpunkte der $2$-dimensionalen Zellen sind. 
}
\end{theorem}

\begin{remark}
\large{Nach Satz 1 ist es klar, dass eine  Einbettungsobstruktionsheorie von 2-dimensionalen Komplexen in $\mathbb{R}^4$, die nur Schnittpunkte von verschiedenen Simplices  ber\"ucksichtigt, nicht vollst\"andig sein kann.
}
\end{remark}

\large{
Dass es $2$-dimensionale Komplexe gibt, die in $\mathbb{R}^4$  nicht eingebettet aber fast-eingebettet werden k\"onnen, wurde von Segal-Skopenkov-Spiez (\cite{SSS}) beobachtet. Aber ihre Beispiele, die sie aus $K_\varphi$ durch die Identifizierung von zwei Punkten $x_1 \sim y_1$ erhalten, sind f\"ur diesen Zweck ein bisschen k\"unstlich, weil die Zellen $x_1 x_2 x_3$ und $y_1 y_2 y_3$ nach dieser Identifizierung  nicht mehr getrennt sind, aber die erhaltenen Komplexe $K_\varphi/( x_1 \sim y_1) $ nicht einfacher eingebettet werden k\"onnen, weil es einfach ist zu zeigen, dass wenn ein solcher Komplex $K_\varphi/(x_1 \sim y_1)$ in $\mathbb{R}^4$ eingebettet werden kann, dann geht dies auch f\"ur $K_\varphi$. 

\subsection*{Beweis} 
\large{
Wir konstruieren eine Immersion $K_\varphi \rightarrow \mathbb{R}^4$, so dass alle Singularit\"aten  Selbstschnittpunkte  der $2$-dimensionalen Zellen $D, x_1x_2x_3$, und $y_1y_2y_3$ sind.

Sei $H$ das Komplement in $K_\varphi$ dieser drei offenen Zellen. Es ist sehr einfach, eine Einbettung $H\rightarrow \mathbb{R}^4$ zu konstruieren, die wir im Folgenden erkl\"aren. Da die $3$-dimensionale Sph\"are $\mathbb{S}^3$ ein Verbund von zwei Kreisen ist, und $\hat{X}$ isomorph zu einem Unterkomplex des Verbundes des Dreiecks $\triangle(x_1, x_2, x_3)$ mit dem Dreieck $\triangle (x_4,x_5,x_6)$ ist, die jeweils hom\"oomorph zu einem Kreis sind, kann der Unterkomplex $\hat{X}$ von $H$, der  aus allen Simplices mit Eckpunkten $x_1,x_2,..., x_6$ aufgebaut ist, in $\mathbb{S}^3$ eingebettet werden.  Diese Einbettung $\hat{X} \rightarrow \mathbb{S}^3$ kann durch einen (abteilungsweise linearen) Selbsthom\"oomorphismus von $\mathbb{S}^3$ so abge\"andert werden, dass das Bild von $\hat{X}$ in einem kleinen Ball $B_X$ enthalten ist, und der Punkt $x_6$ im Rand $\partial B_X$ liegt. \"Ahnlich konstruieren wir  in einem Ball $B_Y \subset \mathbb{S}^3$, der disjunkt von $B_X$ ist,  eine Einbettung eines Unterkomplexes $\hat{Y}$ von $H$, der aus Simplices mit Eckpunkten $y_1,y_2,..., y_6$ aufgebaut ist, so dass der Punkt $y_6$ im Rand $\partial B_Y$ liegt. Dann schieben wir die beiden B\"allen $B_X$ und $B_Y$ aufeinander zu, bis $x_6$ und $y_6$ identifiziert sind, wodurch  wir eine Einbettung $\hat{X}\sqcup\hat{Y}/{x_6 \sim y_6}\rightarrow\mathbb{S}^3 = \partial \mathbb{D}^4$ erhalten. Schlie"slich k\"onnen wir eine Einbettung  $f\colon H \rightarrow \mathbb{R}^4$ erhalten, indem wir zwei Verbunde im $\mathbb{R}^4$ konstruieren, von denen einer ein Verbund des eingebetteten $\hat{X}$ mit einem inneren Punkt in $\mathbb{D}^4$ ist, der das Bild von $x_0$ wird, und der andere ein Verbund des eingebetteten $\hat{Y}$ mit einem inneren Punkt in $\mathbb{D}^4$ ist, der das Bild von $y_0$ wird, so dass die beiden Verbunde sich nicht schneiden.

Es verbleibt jetzt nur noch eine Abbildung von den drei Zellen $D$, $x_1x_2x_3$ und $y_1y_2y_3$ nach $(\mathbb{R}^4 - \text{Int}(\mathbb{D}^4))$ zu konstruieren,
die eine Erweiterung der Einschr\"ankung von $f$ auf deren R\"andern ist, so dass die Bilder dieser drei Zellen paarweise disjunkt sind. Da $(\mathbb{R}^4 - \text{Int}(\mathbb{D}^4))$ und $\mathbb{S}^3 \times [0,\infty)$ hom\"oomorph sind, gibt es f\"ur jede Homotopie $H_t$ in $\mathbb{S}^3$ eine Abbildung $(H_t, t)$ nach $\mathbb{S}^3 \times [0,\infty) \cong (\mathbb{R}^4 - \text{Int}(\mathbb{D}^4))$. Die Abbildung, die noch zu konstruieren ist, ist eine Aneinanderh\"angung von drei Homotopien $F_t$ $G_t$, und $H_t$ in $\mathbb{S}^3$, die wir im Folgenden beschreiben. 
Wir bezeichen durch $\gamma_1$ (bzw. $\gamma_2$)  das Bild des Randes der $2$-dimensionalen Zelle $x_1x_2x_3$ (bzw. $y_1y_2y_3$), und durch $\gamma_3$ eine geschlossene Kurve in $\mathbb{S}^3$, die in einer kleinen Umgebung von $f(\partial D)$ liegt, so dass es eine Homotopie $F_t$ in $\mathbb{S}^3$ zwischen $F_0 = \gamma_3$ und der Anheftungsabbildung $F_1 = f\circ \varphi$ gibt, f\"ur die f\"ur jedes $t \in [0,1)$ die Abbildung $F_t$ eine Einbettung ist. Nach der Arbeit von Milnor (\cite{Milnor}) gibt es eine Verschlingunghomotopie\footnote{d.h. eine Homotopie, w\"ahrend der jede Kurve sich sebst schneiden darf, aber nicht eine andere Kurve schneidet.}  $G_t$ zwischen $\gamma_1, \gamma_2, \gamma_3$ und der trivialen Verschlingung, genau dann wenn die drei Verschlingungszahlen Lk$(\gamma_1,\gamma_2)$, Lk$(\gamma_2,\gamma_3)$,  Lk$(\gamma_3,\gamma_1)$ und die Dreifachverschlingungszahl Lk$(\gamma_1, \gamma_2, \gamma_3)$  gleich Null sind. Die Verschlingungszahl Lk$(\gamma_1,\gamma_2) =0$, weil $\gamma_1$ und $\gamma_2$ in disjunkten B\"allen in $\mathbb{S}^3$ liegen. Die Verschlingungszahlen  Lk$(\gamma_2,\gamma_3)$,  Lk$(\gamma_1,\gamma_3)$ sind gleich Null, wenn $\varphi$ in $[F,F]$ liegt.   Die Dreifachverschlingungszahl von $\gamma_1, \gamma_2, \gamma_3$ ist gleich Null, wenn $\varphi \in F^{(3)} \leq [F,F]$. Die letzte Homotopie $H_t$ schrumpft die triviale Verschlingung zu drei Punkten, so dass $H_t$ f\"ur $t\in [0,1)$ eine Einbettung ist. 
Hiermit haben wir endlich die erw\"unschte Immersion $K_\varphi \rightarrow \mathbb{R}^4$ erhalten, die aus $f$ und den drei Homotopien $F_t, G_t, H_t$ konstruiert wird. Das beendet den Beweis.


} 

\subsection*{Danksagungen}
Ich bin dankbar f\"ur meine Zeit beim Max-Planck-Institut f\"ur Mathematik in Bonn und beim Karlsruher Institut f\"ur Technologie. Ich danke Herr JProf. Dr. Claudio Llosa Isenrich daf\"ur, dass er die vielen Fehler in fr\"uheren Versionen dieses Artikels gefunden und korrigiert hat, und mir Eigenheiten der Deutschen Sprache, z.B. den Unterschied  zwischen ,,endlich`` und ,,schlie"slich``, oder  ,,eigentlich`` und ,,tats\"achlich``, erkl\"art hat. 
\bigskip

\bibliography{Reference}
\bibliographystyle{amsplain}

\end{document}